\documentclass[12pt,a4paper]{article}

\usepackage{amsmath,amssymb,amsthm}

\def\st{\mathchoice{:}{:}{\,:\,}{\,:\,}}

\def\size#1{\lvert#1\rvert}

\def\supp{\operatorname{supp}}
\def\ssm{\smallsetminus}
\def\restrictedto{\mathbin{\upharpoonright}}
\newcommand{\forcestext}[2][{}]{\Vdash_{#1}\text{``}{#2}\text{''}}
\theoremstyle{plain}
\newtheorem{thm}{Theorem}[section]
\newtheorem{lem}[thm]{Lemma}

\newtheorem{prop}[thm]{Proposition}

\newtheorem{q}[thm]{Question}

\theoremstyle{definition}
\newtheorem{defn}[thm]{Definition}
\theoremstyle{remark}

\newtheorem{remark}{Remark}

\begin{document}

\title{Remarks on 
	the preservation of topological covering properties 
	under Cohen forcing}
\author{Masaru Kada\thanks{%
	Supported by 
	Grant-in-Aid for Young Scientists (B) 21740080, MEXT.}
%	\thanks{email: \textsf{kada@mi.s.osakafu-u.ac.jp}}
	}
%\date{February 11, 2010}
\date{}
\maketitle

\begin{abstract}
Iwasa 
investigated 
the preservation of 
various covering properties of topological spaces 
under Cohen forcing. 
By 
%refining 
improving 
the argument in Iwasa's paper, 
we prove that the Rothberger property, 
the Menger property 
and selective screenability 
are also preserved under Cohen forcing 
and forcing with the measure algebra. 
\end{abstract}

\section{Introduction}

Let $(X,\tau)$ be a topological space 
and $\mathbb{P}$ a forcing notion. 
In the forcing extension by $\mathbb{P}$, 
the collection $\tau$ of subsets of $X$ 
is no longer a topology on $X$, 
since there are more infinite subsets of $\tau$, 
whose union may not belong to $\tau$. 
However, $\tau$ is still a base of a topology on $X$. 
So we let $\tau^{\mathbb{P}}$ denote a $\mathbb{P}$-name 
for a topology on $\check{X}$ which is generated by $\tau$, 
and we consider $(\check{X},\tau^{\mathbb{P}})$ 
as a topological space corresponding to $(X,\tau)$ 
in the forcing extension. 

For a property $\Phi$ of a topological space, 
%(such as paracompactness, Lindel\"{o}fness, and so on), 
We say 
a forcing notion $\mathbb{P}$ 
\emph{preserves} 
%a property 
$\Phi$ 
%of a topological space 
if, 
whenever $(X,\tau)$ satisfies $\Phi$, 
we have 
$\forcestext[\mathbb{P}]{
	(\check{X},\tau^{\mathbb{P}})\text{ satisfies }\Phi}$.

Grunberg, Junqueira and Tall \cite{GJT:forcingnormality} 
proved that 
Cohen forcing 
preserves paracompactness. 
%and Lindel\"{o}fness. 
Using their ideas, 
Iwasa \cite{Iwasa:cohen} 
extensively studied 
the preservation of 
various covering properties of topological spaces 
under Cohen extensions, 
and proved that Cohen forcing preserves 
the following properties: 
paracompactness,
subparacompactness,
screenability,
$\sigma$-meta\-compactness,
$\sigma$-para\-Lindel\"{o}fness,
Lindel\"{o}fness 
and 
meta\-Lindel\"{o}fness.
%\begin{enumerate}
%\item paracompactness,
%\item subparacompactness,
%\item screenability,
%\item $\sigma$-meta\-compactness,
%\item $\sigma$-para\-Lindel\"{o}fness,
%\item Lindel\"{o}fness,
%\item metaLindel\"{o}fness.
%\end{enumerate}
It is not so hard to observe, 
though it is not explicitly stated, 
that we can prove the same preservation results 
also for the measure algebra. 

In the present paper, 
we will improve the idea used in Iwasa's paper 
and prove that the following properties are also preserved 
under Cohen forcing and forcing with the measure algebra: 
\begin{enumerate}
\item the Rothberger property,
\item the Menger property,
\item selective screenability.
\end{enumerate}

A space $(X,\tau)$ 
\emph{has the Rothberger property}
if, 
for every sequence $\langle\mathcal{U}_n\st n<\omega\rangle$ 
of open covers of $(X,\tau)$ 
there is 
a sequence $\langle U_n\st n<\omega\rangle$ 
of open sets of $(X,\tau)$ 
such that 
\begin{itemize}
\item for all $n<\omega$, $U_n\in\mathcal{U}_n$, and 
\item $\{U_n\st n<\omega\}$ is an open cover of $(X,\tau)$. 
\end{itemize}
A space $(X,\tau)$ 
\emph{has the Menger property}
if, 
for every sequence $\langle\mathcal{U}_n\st n<\omega\rangle$ 
of open covers of $(X,\tau)$ 
there is 
a sequence $\langle \mathcal{F}_n\st n<\omega\rangle$ 
of 
sets 
%collections 
of open sets of $(X,\tau)$ 
such that 
\begin{itemize}
\item for all $n<\omega$, 
	$\mathcal{F}_n$ is a finite subset of $\mathcal{U}_n$, and 
\item $\bigcup_{n<\omega}\mathcal{F}_n$ is an open cover of $(X,\tau)$. 
\end{itemize}
It is easy to see that 
the Rothberger property 
implies 
the Menger property, 
and the Menger property implies Lindel\"{o}fness. 

For a topological space $(X,\tau)$ 
and two 
%subsets 
sets 
%collections 
$\mathcal{A},\mathcal{B}$ 
%of $\tau$, 
of open sets of $(X,\tau)$, 
we say \emph{$\mathcal{A}$ refines $\mathcal{B}$}, 
or  $\mathcal{A}$ is a \emph{refinement} of $\mathcal{B}$, 
if 
for each $U\in\mathcal{A}$ there is $V\in\mathcal{B}$ with $U\subseteq V$. 
%In this paper, 
We will use this terminology even if $\mathcal{A}$ or $\mathcal{B}$ 
is not a cover of $(X,\tau)$. 

A space $(X,\tau)$ 
is 
\emph{selectively screenable}\footnote{%
	Selective screenability 
	was introduced  
	by Addis and Gresham 
%	in 1978 
	\cite{AG:infdim1} 
	and called property C\null. 
	The term ``selectively screenable'' 
	was coined by Babinkostova in her papers, 
	to avoid the confusion with strong measure zero, 
	which was also called property C in the old days. 
}
if, 
for every sequence $\langle\mathcal{U}_n\st n<\omega\rangle$ 
of open covers of $(X,\tau)$, 
there is 
a sequence $\langle\mathcal{H}_n\st n<\omega\rangle$ 
of 
sets 
%collections 
of open sets of $(X,\tau)$ 
%of subsets of $\tau$ 
such that 
\begin{itemize}
\item for all $n<\omega$, 
	$\mathcal{H}_n$ is pairwise disjoint and refines $\mathcal{U}_n$, and 
\item $\bigcup_{n<\omega}\mathcal{H}_n$ is an open cover of $(X,\tau)$. 
\end{itemize}
%It is clear that 
%a selectively screenable space is screenable, 
%and 
It is easy to see that 
a space with the Rothberger property is selectively screenable. 

For an infinite cardinal $\kappa$, 
let $\mathbb{C}(\kappa)$ denote the Cohen forcing notion 
with the index set $\kappa$ 
(that is, 
$\mathbb{C}(\kappa)=\operatorname{Fn}(\kappa,2)$), 
and $\mathbb{B}(\kappa)$ denote 
the measure algebra on $2^\kappa$.

\begin{remark}
Scheepers and Tall 
proved that, 
	if $(X,\tau)$ is a Lindel\"{o}f space 
	and $\kappa\geq\aleph_1$, 
	then 
	$
	\forcestext[\mathbb{C}(\kappa)]
	{(\check{X},\tau^{\mathbb{P}})\text{ has the Rothberger property}}$ 
	\cite[Theorem~11]{ST:lindelof}. 
Since the Rothberger property implies Lindel\"{o}fness, 
%this means 
their result 
%implies 
means 
%the preservation of the Rothberger property 
%under 
%forcing with $\mathbb{C}(\kappa)$ 
that, 
for $\kappa\geq\aleph_1$,
forcing with $\mathbb{C}(\kappa)$ 
preserves 
the Rothberger property. 
They also proved 
that, 
for any infinite cardinal $\kappa$, 
forcing with $\mathbb{B}(\kappa)$ 
preserves 
the Rothberger property 
%	the preservation of 
%	the Rothberger property 
%	is preserved 
%	under forcing with $\mathbb{B}(\kappa)$ 
%	for any infinite cardinal $\kappa$
	\cite[Theorem~15]{ST:lindelof}. 
%Our result 
Theorem~\ref{thm:preserveroth} 
in the present paper 
%will fill 
fills the ``missing part'', 
that is, 
the preservation of the Rothberger property 
under forcing with $\mathbb{C}(\aleph_0)$. 
It is a still unsolved problem 
whether forcing with $\mathbb{C}(\aleph_0)$ 
converts a Lindel\"{o}f space in the ground model 
into a space with the Rothberger property. 

Scheepers and Tall 
also proved 
that 
the Menger property 
is preserved under forcing with $\mathbb{B}(\kappa)$ 
	for any infinite cardinal $\kappa$ 
\cite[Theorem~42]{ST:lindelof}, 
but it has not been stated that 
the Menger property 
is preserved under Cohen forcing. 
Theorem~\ref{thm:preservemenger} 
in the present paper 
was pointed out by Scheepers. 
\end{remark}

\section{Endowments and approximation of open covers}\label{sec:endow}

The 
combinatorial concept 
%notion
of \emph{endowments}, 
also spelled \emph{$n$-dowments}\footnote{%
%The author of the present paper 
%guesses that 
%it 
%is 
The term ``$n$-dowment'' 
was 
%named, 
coined, 
probably after Alan~Dow, 
by one of the authors of the paper \cite{DTW:normalmoore1} 
other than Dow himself. 
Some other people, including Fleissner, 
have called the same structure a \emph{lynx}. 
}, 
of the Cohen forcing notion $\mathbb{C}(\kappa)$ 
was originally invented by Dow \cite{Dow:remote}, 
and 
used in 
%the papers 
\cite{DTW:normalmoore1} and \cite{GJT:forcingnormality} 
in connection with 
the preservation 
of topological properties 
under Cohen forcing. 

%Using endowments, 
%we can nicely approximate an open cover $\dot{\mathcal{U}}$ 
%of $(\check{X},\tau^{\mathbb{C}(\kappa)})$ 
%by 
%%a countable sequence of open covers 
%an open cover 
%of $(X,\tau)$ 
%in the ground model. 

Although 
endowments are originally defined only 
%with 
for 
the Cohen forcing notion, 
%$\mathbb{C}(\kappa)$, 
here we redefine endowments for a general forcing notion.

\begin{defn}\label{defn:endow}
%Suppose that 
%$\mathbb{P}$ is a forcing notion 
%and 
%$\mathbb{P}$ is decomposed into a countable increasing union, 
%say $\mathbb{P}=\bigcup_{n<\omega}\mathbb{P}_n$ 
%where $\mathbb{P}_{n+1}\supseteq \mathbb{P}_{n}$ for each $n$. 
%We assume that such a decomposition of $\mathbb{P}$ 
%is fixed and clear from the context. 
%
%We say 
A forcing notion 
$\mathbb{P}$ is \emph{endowed} 
if there are
a decomposition of $\mathbb{P}$ 
into an increasing union of length $\omega$, 
say $\mathbb{P}=\bigcup_{n<\omega}P_n$ 
where $P_{n}\subseteq P_{n+1}$ for all $n$, 
and 
a sequence $\langle\mathcal{L}_n\st n<\omega\rangle$ 
of sets with the following properties: 
For each $n<\omega$, 
\begin{enumerate}
\item $\mathcal{L}_n$ is a set of finite antichains in $\mathbb{P}$, 
\item\label{item:maxac} 
	for every maximal antichain $A$ in $\mathbb{P}$, 
	there is 
%	a member $L$ of $\mathcal{L}_n$ 
	$L\in\mathcal{L}_n$ 
	with $L\subseteq A$, and
\item\label{item:compatibility} 
	for any $p\in P_n$ and 
	any $n$ elements $L_0,\ldots,L_{n-1}$ of $\mathcal{L}_n$, 
	there are $q_0,\ldots,q_{n-1}$ such that, 
%	each $q_i$ for $i<n$ belongs to $L_i$, 
	$q_i\in L_i$ for each $i<n$, 
	and the set $\{p,q_0,\ldots,q_{n-1}\}$ 
	has a 
%	common 
	lower bound in $\mathbb{P}$. 
\end{enumerate}
We call a sequence $\langle\mathcal{L}_n\st n<\omega\rangle$ 
which meets the above requirements 
a \emph{sequence of endowments} of $\mathbb{P}$, 
and we will say 
\emph{$\mathbb{P}$ is 
endowed with $\langle\mathcal{L}_n\st n<\omega\rangle$}. 
We call each $\mathcal{L}_n$ an \emph{endowment}, 
an \emph{$n$-dowment}\footnote{%
The letter $n$ in the term ``$n$-dowment'' 
may, \emph{but does not have to}, be considered as a parameter, 
only when one set $\mathcal{L}_n$ 
for a specific natural number $n$ is mentioned 
and the very specific alphabet ``$n$'' (not $k$, $i$, etc.) is chosen 
as a variable.  
Otherwise we should regard the letter $n$ just as a part of the name 
and not for a variable. } 
or an \emph{$n$-th endowment}. 
\end{defn}

%For a fixed infinite cardinal $\kappa$, 
%we decompose $\mathbb{C}(\kappa)$ 
%into the increasing union 
%$\mathbb{C}(\kappa)=\bigcup_{n<\omega}C_n$ 
%where 
%$C_n=\{p\in\mathbb{C}(\kappa)\st \size{p}\leq n\}$. 
Cohen forcing notion 
$\mathbb{C}(\kappa)$ 
is typically decomposed 
into the increasing union 
$\mathbb{C}(\kappa)=\bigcup_{n<\omega}C_n$ 
where 
$C_n=\{p\in\mathbb{C}(\kappa)\st \size{p}\leq n\}$, 
and it is actually endowed with respect to this decomposition. 
The following result is called ``Dow's Lemma'' 
\textup{\cite[Lemma~1.1]{DTW:normalmoore1}}.

\begin{thm}
For any infinite cardinal $\kappa$, 
the Cohen forcing notion $\mathbb{C}(\kappa)$ is endowed. 
\end{thm}

However, the clause (\ref{item:compatibility}) 
in the definition of endowments 
is too strong for our purpose in the present paper. 
So we relax the clause (\ref{item:compatibility}) 
and define the notion of 
\emph{weak endowments}.

\begin{defn}\label{defn:weakendow}
%We say 
A forcing notion 
$\mathbb{P}$ is \emph{weakly endowed} 
if there are 
a decomposition of $\mathbb{P}$ 
into an increasing union of length $\omega$, 
say $\mathbb{P}=\bigcup_{n<\omega}P_n$ 
where $P_{n}\subseteq P_{n+1}$ for all $n$, 
and 
a sequence $\langle\mathcal{L}_n\st n<\omega\rangle$ 
of sets with the following properties: 
For each $n<\omega$, 
\begin{enumerate}
\item $\mathcal{L}_n$ is a set of finite antichains in $\mathbb{P}$, 
\item 
	for every maximal antichain $A$ in $\mathbb{P}$, 
	there is 
%	a member $L$ of $\mathcal{L}_n$ 
	$L\in\mathcal{L}_n$ 
	with $L\subseteq A$, and
\item[($\ref{item:compatibility}'$)] 
	for any $p\in P_n$ and 
	$L\in\mathcal{L}_n$, 
	there is $q\in L$  
	such that $p,q$ are compatible in $\mathbb{P}$. 
\end{enumerate}
We call a sequence $\langle\mathcal{L}_n\st n<\omega\rangle$ 
which meets the above requirements 
a \emph{sequence of weak endowments} of $\mathbb{P}$, 
and we will say 
\emph{$\mathbb{P}$ is 
	weakly endowed with $\langle\mathcal{L}_n\st n<\omega\rangle$}. 
We call each $\mathcal{L}_n$ a \emph{weak endowment} 
%a \emph{weak $n$-dowment}
or an \emph{$n$-th weak endowment}. 
\end{defn}

%Even though the following proposition 
%is an immediate consequence of Dow's Lemma, 
The following proposition is essentially proved in 
a sublemma \cite[Lemma~1.0]{DTW:normalmoore1} for the proof of Dow's Lemma. 
%However, 
%we will present a proof here for self-containedness. 
For self-containedness, 
we will present a proof in Appendix. 

\begin{prop}\label{prop:cohenweakendow}
For any infinite cardinal $\kappa$, 
the Cohen forcing notion $\mathbb{C}(\kappa)$ is weakly endowed. 
\end{prop}

We can see that the measure algebra is also weakly endowed. 

\begin{thm}
For any infinite cardinal $\kappa$, 
the measure algebra $\mathbb{B}(\kappa)$ is weakly endowed. 
\end{thm}

\begin{proof}
Just decompose $\mathbb{B}(\kappa)$ 
into the increasing union 
$\mathbb{B}(\kappa)=\bigcup_{n<\omega}B_n$ 
where 
$B_n=\{p\in\mathbb{B}(\kappa)\st \mu(p)\geq 2^{-n}\}$ 
and, 
for each $n$, 
let 
$\mathcal{L}_n$ be the collection of all finite antichains in 
$\mathbb{B}(\kappa)$ 
whose total measure is greater than $1-2^{-n}$. 
\end{proof}

%Suppose that a forcing notion $\mathbb{P}$ is weakly endowed, 
%and fix a corresponding decomposition $\mathbb{P}=\bigcup_{n<\omega}P_n$ 
%and a sequence 
%$\langle \mathcal{L}_n\st n<\omega\rangle$ of weak endowments. 

Let $(X,\tau)$ be a topological space. 
Using 
%weak 
endowments of $\mathbb{C}(\kappa)$, 
we can nicely approximate an open cover $\dot{\mathcal{U}}$ 
of $(\check{X},\tau^{\mathbb{C}(\kappa)})$ 
by 
%a countable sequence of open covers 
an open cover 
of $(X,\tau)$ 
in the ground model, 
%. 
%We review the idea of approximating an open cover $\dot{\mathcal{U}}$ 
%of $(\check{X},\tau^{\mathbb{P}})$ 
%by 
%an open cover of $(X,\tau)$ in the ground model, 
which was 
the idea 
used 
in \cite{DTW:normalmoore1} and \cite{GJT:forcingnormality}. 
We review this idea, 
in a generalized representation for weakly endowed forcing notions. 

Throughout the rest of the present paper, 
we assume that $\mathbb{P}$ is a weakly endowed forcing notion, 
and fix a corresponding decomposition 
$\mathbb{P}=\bigcup_{n<\omega}P_n$ 
and a sequence $\langle\mathcal{L}_n\st n<\omega\rangle$ 
of weak endowments of $\mathbb{P}$.

Let $\dot{\mathcal{U}}$ be 
a $\mathbb{P}$-name for an open cover of $(\check{X},\tau^{\mathbb{P}})$. 
For each $n<\omega$, 
we will construct 
an open cover 
$\mathcal{V}_n(\dot{\mathcal{U}})$ of $(X,\tau)$, 
which we will call 
the \emph{$n$-th approximation of\/ $\dot{\mathcal{U}}$ 
with respect to $\langle \mathcal{L}_n\st n<\omega\rangle$}. 

Since $\dot{\mathcal{U}}$ is forced to 
be an open cover of $(\check{X},\tau^{\mathbb{P}})$ 
and $\tau$ is a base for $\tau^{\mathbb{P}}$, 
for each $x\in X$ 
we can find a $\mathbb{P}$-name $\dot{W}_x$ for an element of $\tau$ 
such that 
\[
\forcestext[\mathbb{P}]{
	\check{x}\in\dot{W}_x
	\text{ and }
	\exists U\in\dot{\mathcal{U}}\,(\dot{W}_x\subseteq U)
	}. 
\]
For each $x\in X$, 
choose a maximal antichain $A_x$ in $\mathbb{P}$, 
and an open set $W_{x,p}\in\tau$ for each $p\in A_x$, 
so that $x\in W_{x,p}$ and 
$p\forcestext[\mathbb{P}]{
		\dot{W}_x=\check{W}_{x,p}
	}$ 
for $p\in A_x$.

Now we fix $n<\omega$ and define 
$\mathcal{V}_n(\dot{\mathcal{U}})$ in the following way. 
For each $x\in X$, find $L_{x,n}\in \mathcal{L}_n$ 
such that $L_{x,n}\subseteq A_x$, 
and let $V_{x,n}=\bigcap\{W_{x,p}\st p\in L_{x,n}\}$. 
Then $V_{x,n}$ is an open set containing $x$. 
Let $\mathcal{V}_n(\dot{\mathcal{U}})=\{V_{x,n}\st x\in X\}$. 

The following property of $\mathcal{V}_n(\dot{\mathcal{U}})$ 
is easily observed.

\begin{lem}\label{lem:approx}
%Suppose that $\mathbb{P}$ is weakly endowed with 
%$\langle \mathcal{L}_n\st n<\omega\rangle$, 
%with respect to a decomposition $\mathbb{P}=\bigcup_{n<\omega}P_n$. 
For $n<\omega$, 
let $\mathcal{V}_n(\dot{\mathcal{U}})$ 
be the $n$-th approximation of $\dot{\mathcal{U}}$ 
with respect to $\langle \mathcal{L}_n\st n<\omega\rangle$. 
Then for each $n<\omega$, 
for any $V\in\mathcal{V}_n(\dot{\mathcal{U}})$ and $p\in P_n$ 
there is $r\in\mathbb{P}$ 
such that 
$r\leq p$ and 
$r\forcestext[\mathbb{P}]{\,
	\exists U\in\dot{\mathcal{U}}\,(\check{V}\subseteq U )}$.
\end{lem}

\begin{proof}
Fix $n<\omega$, 
$V\in\mathcal{V}_n(\dot{\mathcal{U}})$ 
and $p\in P_n$. 
Find $x\in X$ so that $V=V_{x,n}$ in the construction of 
$\mathcal{V}_n(\dot{\mathcal{U}})$, 
and look at $L_{x,n}$. 
By the property ($\ref{item:compatibility}'$) 
of the $n$-th weak endowment $\mathcal{L}_n$, 
find $q\in L_{x,n}$ and $r\in\mathbb{P}$ so that 
$r\leq p$ and $r\leq q$. 
By the definition of $V_{x,n}$, 
we have $V=V_{x,n}\subseteq W_{x,q}$. 
Since 
$\forcestext[\mathbb{P}]{
	\exists U\in\dot{\mathcal{U}}\,(\dot{W}_x\subseteq U)}$ 
and 
$q\forcestext[\mathbb{P}]{\dot{W}_x=\check{W}_{x,q}}$, 
we have 
$	r\forcestext[\mathbb{P}]{
		\exists U\in\dot{\mathcal{U}}\,(\check{V}\subseteq U)
	}$.
\end{proof}

%Although we defined the notion of weak endowments 
%in a general fashion, 
%%in Section~\ref{sec:endow}, 
%we do not have any example of weakly endowed forcing notions 
%other than Cohen forcing notion and the measure algebra so far. 
%
%\begin{q}
%Are there any further example 
%of an endowed, or weakly endowed, forcing notion? 
%\end{q}

\section{Preservation of covering properties}

Iwasa established 
the following result about the preservation of 
covering properties under Cohen forcing \cite[Corollary~2.6]{Iwasa:cohen}. 
Although he actually dealt only with Cohen forcing, 
the proof works for weakly endowed forcing notions. 

\begin{thm}
The following covering properties are preserved 
under forcing with a weakly endowed forcing notion: 
\begin{enumerate}
\item paracompactness,
\item subparacompactness,
\item screenability,
\item $\sigma$-meta\-compactness,
\item $\sigma$-para\-Lindel\"{o}fness,
\item Lindel\"{o}fness,
\item metaLindel\"{o}fness.
\end{enumerate}
\end{thm}

%In this section, we will prove 
This section 
is devoted to the proof of 
the following preservation theorem. 

\begin{thm}
The following covering properties are preserved 
under forcing with a weakly endowed forcing notion: 
\begin{enumerate}
\item the Rothberger property, 
\item the Menger property, 
\item selective screenability.
\end{enumerate}
\end{thm}

The proof will be 
%done 
worked out 
by improving the idea used in Iwasa's paper. 
The following two lemmata 
are inspired by \cite[Lemma~2.3]{Iwasa:cohen}. 

%Throughout this section, 
%we assume that $\mathbb{P}$ is a weakly endowed forcing notion, 
%and fix a corresponding decomposition 
%$\mathbb{P}=\bigcup_{n<\omega}P_n$ 
%and a sequence $\langle\mathcal{L}_n\st n<\omega\rangle$ 
%of weak endowments of $\mathbb{P}$. 

%The following two lemmata 
%will be useful 
%to prove the preservation 
%of both the Rothberger property 
%and selective screenability. 

%For a topological space $(X,\tau)$ 
%and two subfamilies $\mathcal{A},\mathcal{B}\subseteq\tau$, 
%we say \emph{$\mathcal{A}$ refines $\mathcal{B}$}, 
%or  $\mathcal{A}$ is a \emph{refinement} of $\mathcal{B}$, 
%if 
%for each $U\in\mathcal{A}$ there is $V\in\mathcal{B}$ with $U\subseteq V$. 
%We use this terminology even if $\mathcal{A}$ or $\mathcal{B}$ 
%is not a cover of $(X,\tau)$. 

\begin{lem}\label{lem:refine}
Suppose that $(X,\tau)$ is a topological space, 
$\dot{\mathcal{U}}$ is a $\mathbb{P}$-name 
%such that 
%$\forcestext[\mathbb{P}]{\,
%	\dot{\mathcal{U}}\text{ is an open cover of }
%	(\check{X},\tau^{\mathbb{P}})}$, 
for an open cover of $(\check{X},\tau^{\mathbb{P}})$, 
$n<\omega$, 
and $\mathcal{H}\subseteq\tau$ 
is a refinement of $\mathcal{V}_n(\dot{\mathcal{U}})$. 
%, 
%that is, 
%$\forall H\in\mathcal{H}\,\exists V\in\mathcal{V}_n(\dot{\mathcal{U}})\,
%(H\subseteq V)$. 
Then there is a $\mathbb{P}$-name $\dot{\mathcal{W}}$ 
for a subfamily of $\mathcal{H}$ 
with the following properties.
\begin{enumerate}
\item $\forcestext[\mathbb{P}]{\,
	\dot{\mathcal{W}}\text{ refines }\dot{\mathcal{U}}}$. 
\item for any $p\in P_n$ and $H\in\mathcal{H}$ 
	there is $r\in\mathbb{P}$ 
	such that 
	$r\leq p$ 
	and 
	$r\forcestext[\mathbb{P}]{\check{H}\in\dot{\mathcal{W}}}$. 
\end{enumerate}
\end{lem}

%Note that the family $\mathcal{H}$ in the statement 
%is not necessarily a cover of $(X,\tau)$. 

\begin{proof}
Construct a $\mathbb{P}$-name $\dot{\mathcal{W}}$ 
with the following properties: 
For each $H\in\mathcal{H}$, 
for $p\in\mathbb{P}$, 
\begin{itemize}
\item if 
	$p\forcestext[\mathbb{P}]{
		\exists U\in\dot{\mathcal{U}}\,(\check{H}\subseteq U)}$, 
	then 
	$p\forcestext[\mathbb{P}]{
		\check{H}\in\dot{\mathcal{W}}}$, 
	and 
\item if 
	$\forall r\leq p\,
	(r\not\forcestext[\mathbb{P}]{
		\exists U\in\dot{\mathcal{U}}\,(\check{H}\subseteq U)})$ 
	(equivalently, 
	$p\forcestext[\mathbb{P}]{
		\forall U\in\dot{\mathcal{U}}\,(\check{H}\not\subseteq U)}$), 
	then 
	$p\forcestext[\mathbb{P}]{
		\check{H}\notin\dot{\mathcal{W}}}$. 
\end{itemize}
Such construction of $\dot{\mathcal{W}}$ can be done 
using the ``maximal principle'' 
(see \cite[VII Theorem~8.2]{Ku:set}). 
%It is clear 
Use Lemma~\ref{lem:approx} 
to check 
that this $\dot{\mathcal{W}}$ is as desired. 
\end{proof}

\begin{lem}\label{lem:endowrefine}
Let $\langle\dot{\mathcal{U}}_n\st n<\omega\rangle$ 
be a sequence of $\mathbb{P}$-names 
for open covers of 
a space $(\check{X},\tau^{\mathbb{P}})$, 
and $\langle\mathcal{H}_n\st n<\omega\rangle$ 
a sequence of subsets of $\tau$ with the following properties.
\begin{enumerate}
\item For each $n<\omega$, 
	$\mathcal{H}_n$ refines $\mathcal{V}_n(\dot{\mathcal{U}}_n)$.
\item For all $x\in X$, 
	for infinitely many $n<\omega$ there is $H\in\mathcal{H}_n$ 
	such that $x\in H$. 
\end{enumerate}
Then there is a sequence 
$\langle\dot{\mathcal{W}}_n\st n<\omega\rangle$ 
of $\mathbb{P}$-names with the following properties. 
\begin{enumerate}
\item 
	For each $n<\omega$, 
	$\forcestext[\mathbb{P}]{\,
		\dot{\mathcal{W}}_n\subseteq\check{\mathcal{H}}_n
	\text{ and }
	\dot{\mathcal{W}}_n\text{ refines }\dot{\mathcal{U}}_n}$. 
%\item 
%	For each $n<\omega$, 
%	$\forcestext[\mathbb{P}]{\,
%		\dot{\mathcal{W}}_n\text{ refines }\dot{\mathcal{U}}_n}$. 
\item \label{item:alltoghthercover}
	$\forcestext[\mathbb{P}]{\,
		\bigcup_{n<\omega}\dot{\mathcal{W}}_n
		\text{ covers }(\check{X},\tau^{\mathbb{P}})}$. 
\end{enumerate}
\end{lem}

\begin{proof}
For each $n<\omega$, 
construct $\dot{\mathcal{W}_n}$ 
as in Lemma~\ref{lem:refine} 
from $\mathcal{H}_n$. 
We check 
that 
$\dot{\mathcal{W}_n}$'s meet the requirement (\ref{item:alltoghthercover}). 
Fix $x\in X$ and $p\in\mathbb{P}$. 
Choose $m<\omega$ with $p\in P_m$. 
By the assumption, 
we can choose $n\geq m$ and $H\in\mathcal{H}_n$ with $x\in H$. 
Then, 
by Lemma~\ref{lem:refine}, 
there is $r\in\mathbb{P}$ 
such that 
$r\leq p$ 
and 
$r\forcestext[\mathbb{P}]{\check{H}\in\dot{\mathcal{W}_n}}$. 
This means that every $x\in X$ 
is forced to be covered by 
%$\bigcup_{n<\omega}\dot{\mathcal{W}}_n$. 
the union of $\dot{\mathcal{W}}_n$'s. 
\end{proof}

Now we are going to prove the preservation of the Rothberger property. 
We use the equivalent conditions of the Rothberger property 
shown in the following lemma. 
The equivalence $(1)\Leftrightarrow(2)$ is well-known, 
%(proved by taking common refinements of first $n$ covers), 
and 
$(1)\Leftrightarrow(3)$ is easy.

\begin{lem}\label{lem:rothequiv}
For a topological space $(X,\tau)$, 
the following conditions are equivalent: 
\begin{enumerate}
\item
	(the Rothberger property)
	For every sequence 
	$\langle \mathcal{U}_n\st n<\omega\rangle$ 
	of open covers of $(X,\tau)$, 
	there is a sequence 
	$\langle U_n\st n<\omega\rangle$ 
	of open sets of $(X,\tau)$
	such that 
	\begin{itemize}
	\item for all $n<\omega$, $U_n\in\mathcal{U}_n$, and 
	\item for all $x\in X$ 
		there is $n<\omega$ such that $x\in U_n$. 
	\end{itemize}
\item\label{item:rothinf}
	For every sequence 
	$\langle \mathcal{U}_n\st n<\omega\rangle$ 
	of open covers of $(X,\tau)$, 
	there is a sequence 
	$\langle U_n\st n<\omega\rangle$ 
	of open sets of $(X,\tau)$
	such that 
	\begin{itemize}
	\item for all $n<\omega$, $U_n\in\mathcal{U}_n$, and 
	\item for all $x\in X$ 
		there are infinitely many $n<\omega$ such that $x\in U_n$. 
	\end{itemize}
\item\label{item:rothrefine}
	For every sequence 
	$\langle \mathcal{U}_n\st n<\omega\rangle$ 
	of open covers of $(X,\tau)$, 
	there is a sequence 
	$\langle \mathcal{W}_n\st n<\omega\rangle$ 
	of sets of open sets of $(X,\tau)$
	such that 
	\begin{itemize}
	\item for all $n<\omega$, 
		$\size{\mathcal{W}_n}\leq 1$ 
		(that is, $\mathcal{W}_n$ is either a singleton or $\emptyset$)
%	, 
%	\item for all $n<\omega$, 
		and 
		$\mathcal{W}_n$ refines $\mathcal{U}_n$, and 
	\item for all $x\in X$ 
		there is $n<\omega$ such that, 
		there is $U\in\mathcal{W}_n$ 
		such that $x\in U$. 
	\end{itemize}
\end{enumerate}
\end{lem}

\begin{thm}\label{thm:preserveroth}
Suppose that 
$(X,\tau)$ is a topological space with the Rothberger property 
and $\mathbb{P}$ is a weakly endowed forcing notion. 
Then 
we have 
\[
\forcestext[\mathbb{P}]{
	(\check{X},\tau^{\mathbb{P}})\text{ has the Rothberger property}}. 
\]
\end{thm}

\begin{proof}
Fix a sequence 
	$\langle \dot{\mathcal{U}}_n\st n<\omega\rangle$ 
	of 
	$\mathbb{P}$-names for open covers of $(\check{X},\tau^{\mathbb{P}})$. 
Consider the sequence 
	$
%	\vec{\mathcal{V}}=
	\langle \mathcal{V}_n(\dot{\mathcal{U}}_n)\st n<\omega\rangle$ 
	of open covers of $(X,\tau)$. 
%By applying 
Using 
the condition (\ref{item:rothinf}) in Lemma~\ref{lem:rothequiv}, 
%to $\vec{\mathcal{V}}$, 
we can get a sequence 
	$\langle U_n\st n<\omega\rangle$ 
	of open sets of $(X,\tau)$
	such that 
	\begin{itemize}
	\item for all $n<\omega$, $U_n\in\mathcal{V}_n(\dot{\mathcal{U}}_n)$, 
		and 
	\item for all $x\in X$ 
		there are infinitely many $n<\omega$ such that $x\in U_n$. 
	\end{itemize}
Let $\mathcal{H}_n=\{U_n\}$ for each $n<\omega$. 
Now we apply Lemma~\ref{lem:endowrefine} to 
$\langle\mathcal{H}_n\st n<\omega\rangle$ 
to get a sequence 
$\langle\dot{\mathcal{W}}_n\st n<\omega\rangle$. 
It is straightforward to check that 
$\langle\dot{\mathcal{W}}_n\st n<\omega\rangle$ 
is forced to meet the condition 
(\ref{item:rothrefine}) 
in Lemma~\ref{lem:rothequiv}. 
\end{proof}

Scheepers pointed out that 
the preservation of the Menger property 
under forcing with a weakly endowed forcing notion 
is proved in the same way as 
the proof of 
Theorem~\ref{thm:preserveroth}. 
We will use Lemma~\ref{lem:mengerequiv} 
instead of Lemma~\ref{lem:rothequiv}.

\begin{lem}\label{lem:mengerequiv}
For a topological space $(X,\tau)$, 
the following conditions are equivalent: 
\begin{enumerate}
\item
	(the Menger property)
	For every sequence 
	$\langle \mathcal{U}_n\st n<\omega\rangle$ 
	of open covers of $(X,\tau)$, 
	there is a sequence 
	$\langle \mathcal{H}_n\st n<\omega\rangle$ 
	of sets of open sets of $(X,\tau)$
	such that 
	\begin{itemize}
	\item for all $n<\omega$, 
		$\mathcal{H}_n$ is a finite subset of $\mathcal{U}_n$, and 
	\item for all $x\in X$ 
		there is $n<\omega$ such that, 
		there is $U\in\mathcal{H}_n$ such that $x\in U$. 
	\end{itemize}
\item\label{item:mengerinf}
	For every sequence 
	$\langle \mathcal{U}_n\st n<\omega\rangle$ 
	of open covers of $(X,\tau)$, 
	there is a sequence 
	$\langle \mathcal{H}_n\st n<\omega\rangle$ 
	of sets of open sets of $(X,\tau)$
	such that 
	\begin{itemize}
	\item for all $n<\omega$, 
		$\mathcal{H}_n$ is a finite subset of $\mathcal{U}_n$, and 
	\item for all $x\in X$ 
		there are infinitely many $n<\omega$ such that, 
		there is $U\in\mathcal{H}_n$ such that $x\in U$. 
	\end{itemize}
\item\label{item:mengerrefine}
	For every sequence 
	$\langle \mathcal{U}_n\st n<\omega\rangle$ 
	of open covers of $(X,\tau)$, 
	there is a sequence 
	$\langle \mathcal{W}_n\st n<\omega\rangle$ 
	of sets of open sets of $(X,\tau)$
	such that 
	\begin{itemize}
	\item for all $n<\omega$, 
		$\mathcal{W}_n$ is finite
%	, 
%	\item for all $n<\omega$, 
		and 
%		$\mathcal{W}_n$ 
		refines $\mathcal{U}_n$, and 
	\item for all $x\in X$ 
		there is $n<\omega$ such that, 
		there is $U\in\mathcal{W}_n$ 
		such that $x\in U$. 
	\end{itemize}
\end{enumerate}
\end{lem}

\begin{thm}\label{thm:preservemenger}
Suppose that 
$(X,\tau)$ is a topological space with the Menger property 
and $\mathbb{P}$ is a weakly endowed forcing notion. 
Then 
we have 
\[
\forcestext[\mathbb{P}]{
	(\check{X},\tau^{\mathbb{P}})\text{ has the Menger property}}. 
\]
\end{thm}

We turn to the preservation of selective screenability. 
We use the equivalent conditions of selective screenability 
shown in the following lemma, 
which is easy to check.

\begin{lem}\label{lem:selscrequiv}
For a topological space $(X,\tau)$, 
the following conditions are equivalent: 
\begin{enumerate}
\item\label{item:selscrorg}
	(selective screenability)
	For every sequence 
	$\langle \mathcal{U}_n\st n<\omega\rangle$ 
	of open covers of $(X,\tau)$, 
	there is a sequence 
	$\langle \mathcal{W}_n\st n<\omega\rangle$ 
	of sets of open sets of $(X,\tau)$
	such that 
	\begin{itemize}
	\item for all $n<\omega$, 
		$\mathcal{W}_n$ is 
		pairwise disjoint and refines $\mathcal{U}_n$, and 
	\item for all $x\in X$ 
		there is $n<\omega$ and $U\in\mathcal{W}_n$ such that $x\in U$. 
	\end{itemize}
\item\label{item:selscrinf}
	For every sequence 
	$\langle \mathcal{U}_n\st n<\omega\rangle$ 
	of open covers of $(X,\tau)$, 
	there is a sequence 
	$\langle \mathcal{H}_n\st n<\omega\rangle$ 
	of sets of open sets of $(X,\tau)$
	such that 
	\begin{itemize}
	\item for all $n<\omega$, 
		$\mathcal{H}_n$ is pairwise disjoint 
		and refines $\mathcal{U}_n$, and 
	\item for each $x\in X$, 
		there are infinitely many $n<\omega$ such that, 
		there is $U\in\mathcal{H}_n$ with $x\in U$. 
	\end{itemize}
\end{enumerate}
\end{lem}

%We will call a sequence 
%	$\langle \mathcal{W}_n\st n<\omega\rangle$ 
%in the condition (\ref{item:selscrorg}) 
%a \emph{selective screening} of 
%	$\langle \mathcal{U}_n\st n<\omega\rangle$. 

\begin{thm}
Suppose that 
$(X,\tau)$ is a selectively screenable topological space 
and $\mathbb{P}$ is a weakly endowed forcing notion. 
Then 
we have 
\[
\forcestext[\mathbb{P}]{
	(\check{X},\tau^{\mathbb{P}})\text{ is selectively screenable}}. 
\]
\end{thm}

\begin{proof}
Fix a sequence 
	$\langle \dot{\mathcal{U}}_n\st n<\omega\rangle$ 
	of 
	$\mathbb{P}$-names for open covers of $(\check{X},\tau^{\mathbb{P}})$. 
Consider the sequence 
	$
%	\vec{\mathcal{V}}=
	\langle \mathcal{V}_n(\dot{\mathcal{U}}_n)\st n<\omega\rangle$ 
	of open covers of $(X,\tau)$. 
%By applying 
Using 
the condition (\ref{item:selscrinf}) in Lemma~\ref{lem:selscrequiv}, 
%to $\vec{\mathcal{V}}$, 
we can get a sequence 
	$\langle \mathcal{H}_n\st n<\omega\rangle$ 
	of sets of open sets of $(X,\tau)$
	such that 
	\begin{itemize}
	\item for all $n<\omega$, 
		$\mathcal{H}_n$ is pairwise disjoint and 
		refines $\mathcal{V}_n(\mathcal{U}_n)$, and 
	\item for each $x\in X$, 
		there are infinitely $n<\omega$ such that, 
		there is $U\in\mathcal{H}_n$ with $x\in U$. 
	\end{itemize}
Now we apply Lemma~\ref{lem:endowrefine} to 
$\langle\mathcal{H}_n\st n<\omega\rangle$ 
to get a sequence 
$\langle\dot{\mathcal{W}}_n\st n<\omega\rangle$. 
It is straightforward to check that 
$\langle\dot{\mathcal{W}}_n\st n<\omega\rangle$ 
is forced to 
meet the condition the condition (\ref{item:selscrorg}) 
in Lemma~\ref{lem:selscrequiv}. 
%be a selective screening of 
%$(\check{X},\tau^{\mathbb{P}})$.  
\end{proof}

\section{Question}
Although we defined the notion of weak endowments 
in a general fashion in Section~\ref{sec:endow}, 
%$\mathbb{C}(\kappa)$ 
%and 
%$\mathbb{B}(\kappa)$ 
%are only examples of weakly endowed forcing notions which we have so far. 
we do not have any examples of weakly endowed forcing notions 
other than 
$\mathbb{C}(\kappa)$ 
%Cohen forcing notion 
and 
$\mathbb{B}(\kappa)$ 
%the measure algebra 
so far. 

\begin{q}
Are there any further examples 
of weakly endowed forcing notions? 
\end{q}

\section*{Appendix: Endowing Cohen forcing notions}

Here we present a proof of Proposition~\ref{prop:cohenweakendow}, 
which is based on \cite[Lemma~1.0]{DTW:normalmoore1}. 

\begin{proof}[Proof of Proposition~\ref{prop:cohenweakendow}]
Fix an infinite cardinal $\kappa$ 
and decompose $\mathbb{C}(\kappa)$ 
into an increasing union 
$\mathbb{C}(\kappa)=\bigcup_{n<\omega}C_n$ 
where $C_n=\{p\in\mathbb{C}(\kappa)\st\size{p}\leq n\}$. 

Fix $n<\omega$. 
We will claim the following: 
\begin{quote}
For every maximal antichain $A$ in $\mathbb{C}(\kappa)$ 
there is a finite subset $L$ of $A$ such that, 
for every $p\in C_n$ there is $q\in L$ which is compatible with $p$. 
\end{quote}
Then we gather all $L$'s, each corresponding to a maximal antichain, 
and it makes up the $n$-th weak endowment $\mathcal{L}_n$. 

For $p\in\mathbb{C}(\kappa)$, 
$\supp(p)$ denotes 
the domain of $p$ as a partial function from $\kappa$ to $2$, 
and for $F\subseteq\mathbb{C}(\kappa)$, 
let $\supp(F)=\bigcup\{\supp(p)\st p\in F\}$. 

Fix a maximal antichain $A$ in $\mathbb{C}(\kappa)$. 
Pick any $a\in A$. 
Let $E_0=\{a\}$ and $D_0=\supp(E_0)$. 
For each $p\in\mathbb{C}(\kappa)$ with $\supp(p)\subseteq D_0$, 
choose exactly one $a_p\in A$ which is compatible with $p$. 
Let $E_1$ be the collection of such $a_p$'s 
and $D_1=\supp(E_0\cup E_1)$. 
Similarly, 
for $i\leq n$, 
obtain $E_i$ and $D_i$ from $D_{i-1}$. 
Note that $D_i\supseteq D_{i-1}$ holds for each $i$. 
Finally we let $L=\bigcup_{i\leq n}E_i$. 

We check that this $L$ works. 
Fix an arbitrary $p\in P_n$. 
Consider $n+1$ disjoint subsets 
$D_0, D_1\ssm D_0, \ldots, D_{n}\ssm D_{n-1}$ of $D_{n}$. 
Since $\size{\supp(p)}\leq n$, 
$\supp(p)$ must be disjoint from some of those pieces, 
say $D_{i}\ssm D_{i-1}$ (let $D_{-1}=\emptyset$ for convention). 
Then there is $q\in E_i$ which is compatible with 
$p\restrictedto D_{i-1}$. 
Since $\supp(q)\subseteq D_i$ 
and $\supp(p)\cap(D_i\ssm D_{i-1})=\emptyset$, 
$q$ is compatible with $p$. 
\end{proof}

\section*{Acknowledgement}

I thank 
Marion~Scheepers
and 
Franklin~D.~Tall 
for many helpful comments 
on a provisional version of this article. 

%\bibliographystyle{amsplain}
%\bibliography{kada}

\providecommand{\bysame}{\leavevmode\hbox to3em{\hrulefill}\thinspace}
\providecommand{\MR}{\relax\ifhmode\unskip\space\fi MR }
\providecommand{\MRhref}[2]{%
  \href{http://www.ams.org/mathscinet-getitem?mr=#1}{#2}
}
\providecommand{\href}[2]{#2}

\bigskip

\begin{flushleft}
\begin{small}
Masaru Kada\\
Graduate School of Science\\
Osaka Prefecture University\\
1--1 Gakuen-cho, Naka-ku, Sakai, Osaka 599--8531 JAPAN\\
email: \texttt{kada@mi.s.osakafu-u.ac.jp}
\end{small}
\end{flushleft}

\end{document}